\documentclass[a4paper]{article}

\usepackage{amsthm}
\usepackage{amsmath}
\usepackage{booktabs}
\usepackage{color}
\usepackage{enumitem}
\usepackage[hang,flushmargin]{footmisc}
\usepackage[left=25mm,right=25mm,top=25mm,bottom=25mm]{geometry}
\usepackage[
	pdftitle={Quantum curves for the enumeration of ribbon graphs and hypermaps},
	pdfauthor={Norman Do and David Manescu},
	ocgcolorlinks,
	linkcolor=linkblue,
	citecolor=linkred,
	urlcolor=linkred]
{hyperref}
\usepackage{mathpazo}
\usepackage{mathtools}
\usepackage{microtype}
\usepackage{sectsty}
\usepackage{setspace}
\usepackage{titlesec}
	\titlespacing{\section}{0pt}{12pt}{0pt}
	\titlespacing{\subsection}{0pt}{6pt}{0pt}

\theoremstyle{plain}
\newtheorem{theorem}{Theorem}
\newtheorem{proposition}[theorem]{Proposition}

\newtheorem{lemma}[theorem]{Lemma}
\newtheorem{conjecture}[theorem]{Conjecture}
\newtheorem*{conjecture*}{Conjecture}

\theoremstyle{definition}
\newtheorem{definition}[theorem]{Definition}
\newtheorem{example}[theorem]{Example}

\theoremstyle{remark}
\newtheorem{remark}[theorem]{Remark}

\definecolor{linkred}{rgb}{0.75,0,0}
\definecolor{linkblue}{rgb}{0,0,1}

\newcommand\blfootnote[1]{%
  \begingroup
  \renewcommand\thefootnote{}\footnote{#1}%
  \addtocounter{footnote}{-1}%
  \endgroup
}

\sectionfont{\large}
\subsectionfont{\normalsize}

\setlist{nolistsep}

\newcommand {\dd}{\mathrm{d}}
\newcommand {\h}{\hbar}

\makeatletter
\newcommand{\dashrule}[1][black]{
  \color{#1}\rule[\dimexpr.5ex-.2pt]{4pt}{.4pt}\xleaders\hbox{\rule{4pt}{0pt}\rule[\dimexpr.5ex-.2pt]{4pt}{.4pt}}\hfill\kern0pt%
}
\newcommand{\rulecolor}[1]{
  \def\CT@arc@{\color{#1}}
}
\makeatother

\begin{document}

{\large \bfseries Quantum curves for the enumeration of ribbon graphs and hypermaps}

{\bfseries Norman Do and David Manescu}

{\em Abstract.} The topological recursion of Eynard and Orantin governs a variety of problems in enumerative geometry and mathematical physics. The recursion uses the data of a spectral curve to define an infinite family of multidifferentials. It has been conjectured that, under certain conditions, the spectral curve possesses a non-commutative quantisation whose associated differential operator annihilates the partition function for the spectral curve. In this paper, we determine the quantum curves and partition functions for an infinite sequence of enumerative problems involving generalisations of ribbon graphs known as hypermaps. These results give rise to an explicit conjecture relating hypermap enumeration to the topological recursion and we provide evidence to support this conjecture.
\blfootnote{{\em 2010 Mathematics Subject Classification:} 05A15; 14N10; 81S10. \\
{\em Date:} 25 December 2013 \\ The first author was supported by the Australian Research Council grant DE130100650. The second author was supported by a Monash University Vacation Research Scholarship. \\ The first author would like to thank Oliver Leigh and Paul Norbury for valuable discussions.}

~

\hrule

\section{Introduction} \label{introduction}

The topological recursion of Eynard and Orantin was inspired by the loop equations from the theory of matrix models~\cite{eyn-ora07a}. It is known to govern a variety of problems in enumerative geometry and mathematical physics, including intersection theory on moduli spaces of curves~\cite{eyn-ora07a}, Weil--Petersson volumes of moduli spaces of hyperbolic surfaces~\cite{eyn-ora07b}, enumeration of ribbon graphs~\cite{nor13, dum-mul-saf-sor}, stationary Gromov--Witten theory of $\mathbb{P}^1$~\cite{nor-sco, dun-ora-sha-spi}, simple Hurwitz numbers and their generalisations~\cite{bou-mar, eyn-mul-saf, do-lei-nor, bou-her-liu-mul}, and Gromov--Witten theory of toric Calabi--Yau threefolds~\cite{bou-kle-mar-pas, eyn-ora12}. There are also conjectural relations to spin Hurwitz numbers~\cite{mul-sha-spi} and quantum invariants of knots~\cite{dij-fuj-man, bor-eyn}. Given these myriad applications, one would like to determine the scope of universality of the topological recursion, as well as the commonality among problems governed by it.

The topological recursion takes as input the data of a spectral curve --- essentially, a compact Riemann surface $C$ endowed with two meromorphic functions $x$ and $y$. The output is an infinite family of meromorphic multidifferentials $\omega_{g,n}$ on $C^n$ for integers $g \geq 0$ and $n \geq 1$. The solutions to problems from enumerative geometry arise as the coefficients of certain series expansions of these multidifferentials. In this paper, we will be concerned with the rational spectral curves given by the parametrisation
\[
x(z) = z^{a-1} + \frac{1}{z} \qquad \text{and} \qquad y(z) = z, \qquad \text{for a positive integer $a$}.
\]

It was posited by Gukov and Su{\l}kowski that spectral curves $A(x, y) = 0$ satisfying a certain K-theoretic criterion may be quantised to produce a non-commutative curve $\widehat{A}(\widehat{x}, \widehat{y})$. One can interpret $\widehat{A}(\widehat{x}, \widehat{y})$ as a differential operator via $\widehat{x} = x$ and $\widehat{y} = -\h \frac{\partial}{\partial x}$, and it is natural to consider the following Schr\"{o}dinger-like equation~\cite{guk-sul}.
\[
\widehat{A}(\widehat{x}, \widehat{y}) \, Z(x, \h) = 0
\]

Gukov and Su{\l}kowski conjecture that the solution $Z(x, \h)$ is a naturally defined partition function whose perturbative expansion can be calculated from the spectral curve via the topological recursion. Conversely, they suggest that the Schr\"{o}dinger-like equation may be used to recover the quantisation $\widehat{A}(\widehat{x}, \widehat{y})$.

The quantum curve has been rigorously shown to exist in the sense of Gukov and Su{\l}kowski for several problems, including intersection theory on moduli spaces of curves~\cite{zho12a}, enumeration of ribbon graphs~\cite{mul-sul}, simple Hurwitz numbers and their generalisations~\cite{mul-sha-spi}, and open string invariants for $\mathbb{C}^3$ and the resolved conifold~\cite{zho12b}.

In this paper, we consider the following enumerative geometry problem for a fixed positive integer $a$. Let $M_{g,n}^{[a]}(b_1, b_2, \ldots, b_n)$ denote the weighted count of connected genus $g$ branched covers of  marked Riemann surfaces $f: (S; p_1, p_2, \ldots, p_n) \to (\mathbb{P}^1; \infty)$ such that
\begin{itemize}
\item $f$ is unramified over $\mathbb{P}^1 \setminus \{0, 1, \infty\}$;
\item each point in $f^{-1}(1)$ has ramification order $a$; and
\item the preimage divisor $f^{-1}(\infty)$ equals $b_1p_1 + b_2p_2 + \cdots + b_np_n$. 
\end{itemize}
By considering the monodromy of the branched covers, $M_{g,n}^{[a]}(b_1, b_2, \ldots, b_n)$ can be interpreted as a weighted count of certain factorisations in symmetric groups. Such factorisations in turn correspond to certain decorated cell decompositions of a genus $g$ surface, which we refer to as {\em $a$-hypermaps}. These are natural generalisations of ribbon graphs, which are recovered in the case $a = 2$.

One can define generating functions for the enumeration of $a$-hypermaps known as {\em free energies}.
\[
F_{g,n}^{[a]}(x_1, x_2, \ldots, x_n) = \sum_{b_1, b_2, \ldots, b_n = 1}^\infty M_{g,n}^{[a]}(b_1, b_2, \ldots, b_n) \, x_1^{-b_1} x_2^{-b_2} \cdots x_n^{-b_n}
\]
A subtle exception to the equation above is required in the case $(g,n) = (0,1)$.
\[
F_{0,1}^{[a]}(x_1) = -\log x_1 + \sum_{b_1 = 1}^\infty M_{0,1}^{[a]}(b_1) \, x_1^{-b_1}
\]
The {\em partition function} is defined from the free energies in the following way.
\[
Z^{[a]}(x, \h) = \exp \left[ \sum_{g=0}^\infty \sum_{n=1}^\infty \frac{\h^{2g-2+n}}{n!} \, F_{g,n}(x, x, \ldots, x) \right]
\]

The main theorem of this paper is the following.

\begin{theorem} \label{th:qcurve}
The partition function $Z^{[a]}(x, \h)$ satisfies the following Schr\"{o}dinger-like equation for each positive integer $a$, where $\widehat{x} = x$ and $\widehat{y} = -\h \frac{\partial}{\partial x}$.
\[
\left[ \widehat{y}^a - \widehat{x} \widehat{y} + 1 \right] Z^{[a]}(x, \h) = 0
\]
\end{theorem}

This result implies that the quantum curve for the enumeration of $a$-hypermaps is given by the equation $\widehat{A}(\widehat{x},\widehat{y}) = \widehat{y}^a - \widehat{x} \widehat{y} + 1$. This adds an infinite sequence of examples to the growing list of known quantum curves. In the case $a = 2$, we recover the quantum curve for the enumeration of ribbon graphs, which was rigorously determined by Mulase and Su{\l}kowski~\cite{mul-sul}. Our proof bypasses the calculations in their paper by interpreting the Schr\"{o}dinger-like equation as a direct combinatorial statement concerning disconnected unlabelled hypermaps.

The corresponding spectral curve for the enumeration of $a$-hypermaps should arise as the semi-classical limit $A(x, y) = y^a - xy + 1 = 0$. This observation immediately suggests the following conjecture.

\begin{conjecture} \label{con:hypermaps}
For a fixed positive integer $a$, the topological recursion applied to the rational spectral curve $x(z) = z^{a-1} + \frac{1}{z}$ and $y(z) = z$ produces correlation differentials whose expansions at $x_i = \infty$ satisfy
\[
\omega_{g,n} = \sum_{b_1, b_2, \ldots, b_n = 1}^\infty M_{g,n}^{[a]}(b_1, b_2, \ldots, b_n) \, \prod_{i=1}^n \frac{b_i}{x_i^{b_i+1}} \, \dd x_i \qquad \qquad \text{for } 2g-2+n>0.
\]
\end{conjecture}

The conjecture is known to be true in the case $a=2$, which follows from results on the enumeration of lattice points in moduli spaces of curves~\cite{nor13, dum-mul-saf-sor}. There is strong numerical evidence to support the conjecture, in the form of low genus calculations.

The structure of the paper is as follows.
\begin{itemize}
\item In Section~\ref{sec:toprec}, we define the topological recursion of Eynard and Orantin as well as the notion of a quantum curve, in the sense of Gukov and Su{\l}kowski.
\item In Section~\ref{hypermaps}, we define $M_{g,n}^{[a]}(b_1, b_2, \ldots, b_n)$ to be the weighted count of certain genus $g$ branched covers of $\mathbb{P}^1$. We show that these branched covers are in one-to-one correspondence with hypermaps, which are natural generalisations of ribbon graphs. We conclude the section by describing a graphical representation of hypermaps.
\item In Section~\ref{main-proof}, we provide a rigorous interpretation of Theorem~\ref{th:qcurve}, which is required since the definition of the partition function does not yield a convergent power series. We present a combinatorial proof of the theorem, which also gives rise to a concise explicit formula for the partition function.
\item In Section~\ref{conjecture}, we use Theorem~\ref{th:qcurve} to motivate Conjecture~\ref{con:hypermaps} and provide supporting evidence as well as applications.
\end{itemize}

\section{Topological recursion and quantum curves} \label{sec:toprec}

\subsection{Topological recursion} \label{subsec:toprec}

The topological recursion of Eynard and Orantin was inspired by the theory of matrix models~\cite{eyn-ora07a}. It formalises and generalises the loop equations, which are used to calculate perturbative expansions of matrix model correlation functions. The topological recursion uses the data of a spectral curve $C$ to define a family of meromorphic multidifferentials $\omega_{g,n}$ on $C^n$, for integers $g \geq 0$ and $n \geq 1$. In other words, $\omega_{g,n}$ is a meromorphic section of the line bundle $\pi_1^*(T^*C) \otimes \pi_2^*(T^*C) \otimes \cdots \otimes \pi_n^*(T^*C)$ on the Cartesian product $C^n$, where $\pi_i: C^n \to C$ denotes the projection onto the $i$th factor.

{\bf Input.} The input to the recursion is a {\em spectral curve}, which consists of a compact Riemann surface $C$ endowed with two meromorphic functions $x$ and $y$, as well as a choice of a symplectic basis of $H_1(C, \mathbb{Z})$. We require the zeroes of $\mathrm{d}x$ to be simple and distinct from the zeroes of $\mathrm{d}y$. Extensions of the topological recursion to more general spectral curves have appeared in the literature, although they are necessarily more involved than our purposes demand~\cite{bou-eyn}.

{\bf Base cases.} The base cases for the recursion are
\[
\omega_{0,1}(z_1) = y(z_1)~\dd x(z_1) \qquad \text{and} \qquad \omega_{0,2}(z_1, z_2) = B(z_1, z_2).
\]
Here, $B(z_1, z_2)$ is the unique meromorphic bidifferential on $C \times C$ that
\begin{itemize}
\item is symmetric: $B(z_1, z_2) = B(z_2, z_1)$;
\item is normalised on the $A$-cycles of $H_1(C, \mathbb{Z})$: $\oint_{A_i} B(z_1, \,\cdot\,) = 0$; and
\item has double poles without residue along the diagonal $z_1 = z_2$ but is holomorphic away from the diagonal: $B(z_1, z_2) = \frac{\dd z_1~\dd z_2}{(z_1-z_2)^2} + \text{holomorphic}$.
\end{itemize}
The bidifferential $B(z_1, z_2)$ is a natural construction that is sometimes referred to as the fundamental normalised differential of the second kind on $C$.

{\bf Output.} Recursively define the multidifferentials $\omega_{g,n}$ by the following equation, where $S = \{1, 2, \ldots, n\}$ and $z_I = (z_{i_1}, z_{i_2}, \ldots, z_{i_m})$ for $I = \{i_1, i_2, \ldots, i_m\}$.
\[
\omega_{g,n+1}(z_0, z_S) = \sum_{\alpha} \mathop{\text{Res}}_{z=\alpha} \, K(z_0, z) \left[ \omega_{g-1,n+2}(z, \overline{z}, z_S) + \mathop{\sum_{g_1+g_2=g}^\circ}_{I \sqcup J = S} \omega_{g_1, |I|+1}(z, z_I) \, \omega_{g_2, |J|+1}(\overline{z}, z_J) \right]
\]
Here, the outer summation is over the zeroes $\alpha$ of $\dd x$. Since the zeroes are assumed to be simple, there exists a unique non-identity meromorphic function $z \mapsto \overline{z}$ defined on a neighbourhood of $\alpha \in C$ such that $x(\overline{z}) = x(z)$. The symbol $\circ$ over the inner summation denotes the fact that we exclude any terms that involve $\omega_{0,1}$. The kernel $K$ appearing in the residue is defined by the following equation.
\[
K(z_0, z) = \frac{\int_z^{\overline{z}} \omega_{0,2}(z_0, \,\cdot\, )}{2 \left[ y(z) - y(\overline{z}) \right] \dd x(z)}
\]

The multidifferentials $\omega_{g,n}$ have been referred to in the literature as both Eynard--Orantin invariants and correlation functions. Since they are neither functions nor invariant under prescribed transformations, we will use the term {\em correlation differentials}.

The topological recursion has found wide applicability beyond the realm of matrix models, whence it was first conceived. It is now known to govern a variety of problems in enumerative geometry and mathematical physics, with conjectural relations to many more. The following table lists some of these problems and their associated spectral curves, with those yet to be rigorously proven below the dashed line.

\rulecolor{linkblue}
\begin{figure}[ht!]
\centering
\begin{tabular}{@{}p{105mm}@{}ll@{}} \toprule
PROBLEM & \multicolumn{2}{@{}l@{}}{SPECTRAL CURVE} \\ \midrule
intersection theory on moduli spaces of curves~\cite{eyn-ora07a} & $x(z) = z^2$ & $y(z) = z$ \\
enumeration of ribbon graphs~\cite{nor13, dum-mul-saf-sor} & $x(z) = z + \frac{1}{z}$ & $y(z) = z$ \\
Weil--Petersson volumes of moduli spaces~\cite{eyn-ora07b} & $x(z) = z^2$ & $y(z) = \frac{\sin(2\pi z)}{2\pi}$ \\ 
stationary Gromov--Witten theory of $\mathbb{P}^1$~\cite{nor-sco, dun-ora-sha-spi} & $x(z) = z + \frac{1}{z}$ & $y(z) = \log(z)$ \\
simple and orbifold Hurwitz numbers~\cite{bor-eyn-mul-saf, eyn-mul-saf, do-lei-nor, bou-her-liu-mul} & $x(z) = z \exp(-z^a)$ & $y(z) = z^a$ \\
Gromov--Witten theory of toric Calabi--Yau threefolds~\cite{eyn-ora12} \hspace{10pt} & \multicolumn{2}{@{}l@{}}{mirror curves} \\
\multicolumn{3}{@{}c@{}}{\makebox[\linewidth]{\dashrule[linkblue]}} \\
spin Hurwitz numbers~\cite{mul-sha-spi} & $x(z) = z \exp(-z^r)$ & $y(z) = z$ \\
asymptotics of coloured Jones polynomials of knots~\cite{dij-fuj-man, bor-eyn} & \multicolumn{2}{@{}l@{}}{$A$-polynomials} \\ \bottomrule
\end{tabular}
\end{figure}

\subsection{Quantum curves} \label{subsec:qcurves}

Following the work of Gukov and Su{\l}kowski, we use the correlation differentials produced by the topological recursion to define {\em free energies}~\cite{guk-sul}.
\[
F_{g,n}(x_1, x_2, \ldots, x_n) = \int_p^{x_1} \!\! \int_p^{x_2} \! \cdots \! \int_p^{x_n} \omega_{g,n}(x_1, x_2, \ldots, x_n)
\]
We choose a point $p$ on the spectral curve such that $x(p) = \infty$ as the base point for each of the nested integrals~\cite{bor-eyn}. The free energies are in turn used to define a natural {\em partition function}.
\[
Z(x, \h) = \exp \left[ \sum_{g=0}^\infty \sum_{n=1}^\infty \frac{\h^{2g-2+n}}{n!} F_{g,n}(x, x, \ldots, x) \right]
\]
Note that the definition provided here applies only to genus zero spectral curves. In the case of higher genus, it has been proposed in the physics literature that non-perturbative correction terms involving derivatives of theta functions associated to the spectral curve are required~\cite{bor-eyn}.

Given a spectral curve in the form $A(x, y) = 0$, one can ask whether there exists a quantisation $\widehat{A}(\widehat{x}, \widehat{y})$. This quantum curve is non-commutative in the sense that $\widehat{x}$ and $\widehat{y}$ satisfy the commutation relation $[\widehat{x}, \widehat{y}] = \h$. Thus, it is natural for the multiplication operator $\widehat{x} = x$ and the differentiation operator $\widehat{y} = - \h \frac{\partial}{\partial x}$ to be chosen as the polarisation. We call $\widehat{A}(\widehat{x}, \widehat{y})$ a {\em quantum curve} if we recover the spectral curve $A(x, y) = 0$ in the semi-classical limit $\h \to 0$ and if the following Schr\"{o}dinger-like equation is satisfied.
\[
\widehat{A}(\widehat{x}, \widehat{y}) \, Z(x, \h) = 0
\]

Gukov and Su{\l}kowski posit the existence of a quantum curve for a spectral curve $C$ if and only if a certain $K$-theoretic condition is satisfied --- namely, that the tame symbol $\{x, y\} \in K_2(\mathbb{C}(C))$ is a torsion class~\cite{guk-sul}. Note that this condition is automatically satisfied whenever the spectral curve has genus zero. They furthermore combine the calculation of the partition function via the topological recursion along with the Schr\"{o}dinger-like equation in order to solve for $\widehat{A}$ order by order in powers of $\h$.
\[
\widehat{A} = \widehat{A}_0 + \h \widehat{A}_1 + \h^2 \widehat{A}_2 + \cdots
\]
The paper of Gukov and Su{\l}kowski demonstrates the efficacy of this quantisation process by calculating the first few terms of $\widehat{A}$ and using these to predict the form of the quantum curve in several cases of geometric interest~\cite{guk-sul}. This approach amounts to performing quantisation by travelling {\em the long way around} the following schematic diagram.

\begin{center}
\includegraphics{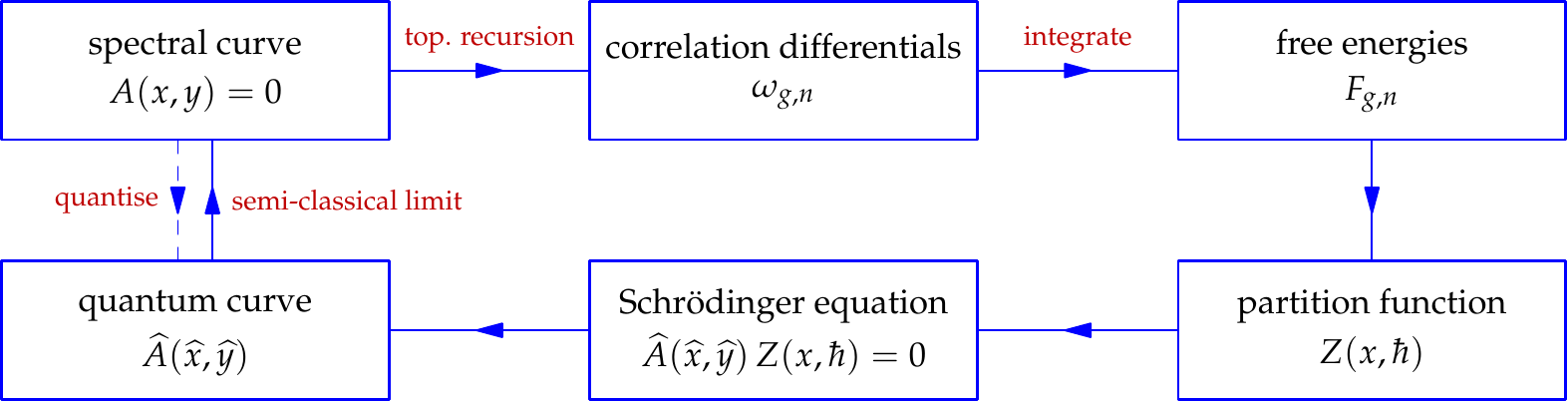}
\end{center}

The quantum curve has been rigorously established for several problems and their associated spectral curves, including intersection theory on moduli spaces of curves~\cite{zho12a}, enumeration of ribbon graphs~\cite{mul-sul}, stationary Gromov--Witten theory of $\mathbb{P}^1$~\cite{dun-mul-nor-pop-sha}, simple Hurwitz numbers and their generalisations~\cite{mul-sha-spi}, and open string invariants for $\mathbb{C}^3$ and the resolved conifold~\cite{zho12b}. The quantum curve for the $A$-polynomial of a knot should recover the $q$-difference operator that appears in the AJ conjecture of Garoufalidis and Le that relates the $A$-polynomial to the coloured Jones polynomials~\cite{gar-le, dij-fuj-man, bor-eyn}.

\section{Ribbon graphs and hypermaps} \label{hypermaps}

\subsection{Counting branched covers of \texorpdfstring{$\mathbb{P}^1$}{P\^{}1}}

In this paper, we consider the following enumeration of branched covers of $\mathbb{P}^1$.

\begin{definition} \label{def:enumeration}
For a fixed positive integer $a$, let $M_{g,n}^{[a]}(b_1, b_2, \ldots, b_n)$ denote the weighted count of connected genus $g$ branched covers of  marked Riemann surfaces $f: (S; p_1, p_2, \ldots, p_n) \to (\mathbb{P}^1; \infty)$ such that
\begin{itemize}
\item $f$ is unramified over $\mathbb{P}^1 \setminus \{0, 1, \infty\}$;
\item each point in $f^{-1}(1)$ has ramification order $a$; and
\item the preimage divisor $f^{-1}(\infty)$ equals $b_1p_1 + b_2p_2 + \cdots + b_np_n$. 
\end{itemize}
\end{definition}

Two branched covers $S_1 \to \mathbb{P}^1$ and $S_2 \to \mathbb{P}^1$ are considered equivalent if there exists an isomorphism $S_1 \to S_2$ that preserves the marked points and covers the identity on $\mathbb{P}^1$. As usual, the weight attached to a branched cover is equal to the reciprocal of its number of automorphisms.

The Riemann existence theorem allows us to represent such a branched cover by its monodromy over $0, 1, \infty \in \mathbb{P}^1$, which motivates the following definition.

\begin{definition} \label{def:hypermap}
For a fixed positive integer $a$, an {\em $a$-hypermap} of type $(g,n)$ is a triple $(\sigma_0, \sigma_1, \sigma_2)$ of permutations in the symmetric group $S_d$ such that
\begin{itemize}
\item $\sigma_0 \sigma_1 \sigma_2 = \text{id}$;
\item $\sigma_2$ consists of $n$ disjoint cycles;
\item $\sigma_1$ has cycle type $(a, a, \ldots, a)$; and
\item $\sigma_0$ has $v$ cycles, where $v = 2 - 2g - n + (a-1) \frac{d}{a}$.
\end{itemize}
\end{definition}

Note that the last condition is simply the constraint implied by the Riemann--Hurwitz formula to ensure that the corresponding branched cover has genus $g$. We say that a hypermap is {\em connected} if the permutations $ \sigma_0, \sigma_1, \sigma_2$ generate a transitive subgroup of $S_d$, which is equivalent to the fact that the corresponding branched cover is connected. We call a hypermap {\em labelled} if the disjoint cycles of $\sigma_2$ are labelled from 1 up to $n$.

Two hypermaps $(\sigma_0, \sigma_1, \sigma_2)$ and $(\tau_0, \tau_1, \tau_2)$ are considered equivalent if one can simultaneously conjugate the $\sigma_i$ to obtain the $\tau_i$ for $i = 0, 1, 2$. If the hypermaps are labelled, then we also impose the condition that the conjugation must preserve the labels. An automorphism of a hypermap is, of course, an isomorphism from a hypermap to itself.

By construction, we may now interpret $M_{g,n}^{[a]}(b_1, b_2, \ldots, b_n)$ as the following weighted count of hypermaps via the Riemann existence theorem.

\begin{proposition}
The number $M_{g,n}^{[a]}(b_1, b_2, \ldots, b_n)$ is equal to the weighted count of connected, labelled $a$-hypermaps $(\sigma_0, \sigma_1, \sigma_2)$ of type $(g,n)$, where the cycle of $\sigma_2$ labelled $i$ has length $b_i$ for $i = 1, 2, \ldots, n$. As usual, the weight attached to a hypermap is equal to the reciprocal of its number of automorphisms.
\end{proposition}

In the case $a = 2$, the notion of an $a$-hypermap reduces to the usual definition of a ribbon graph, which encodes the combinatorics of a cell decomposition of an oriented surface. Ribbon graphs play a prominent role in the study of matrix models and moduli spaces of curves~\cite{lan-zvo}. The ribbon graph case of the enumerative problem considered in Definition~\ref{def:enumeration} has previously appeared in the literature, where it was interpreted as an enumeration of {\em dessins d'enfant}~\cite{dum-mul-saf-sor}.

\subsection{A graphical representation of hypermaps}

As in the case of ribbon graphs, it is useful to represent hypermaps graphically. An $a$-hypermap of type $(g,n)$ corresponds to a bicoloured cell decomposition of a genus $g$ oriented surface. Using the notation of Definition~\ref{def:hypermap}, there are $n$ white faces corresponding to the cycles of $\sigma_2$, where the perimeter of a face is equal to the length of the corresponding cycle. There are also $\frac{d}{a}$ black faces corresponding to the cycles of $\sigma_1$, where the perimeter of each face is equal to $a$. The cell decomposition possesses a bipartite structure in the sense that each edge is incident to one white face and one black face. We informally think of the black faces as $a$-sided edges and refer to them as {\em hyperedges}. Furthermore, we refer to a pair of adjacent edges of a hyperedge as an {\em angle}.

The permutations $\sigma_0, \sigma_1, \sigma_2$ act on the set of $d$ hyperedge angles in the following way. The permutation $\sigma_0$ rotates angles that are adjacent to a common vertex, the permutation $\sigma_1$ rotates angles within a hyperedge, and the permutation $\sigma_2$ rotates angles around faces of the hypermap. These rotations are all performed in an anticlockwise manner with respect to the orientation of the underlying surface. Note that the local geometry of the cell decomposition automatically imposes the condition $\sigma_0 \sigma_1 \sigma_2 = \text{id}$.

\begin{center}
\includegraphics{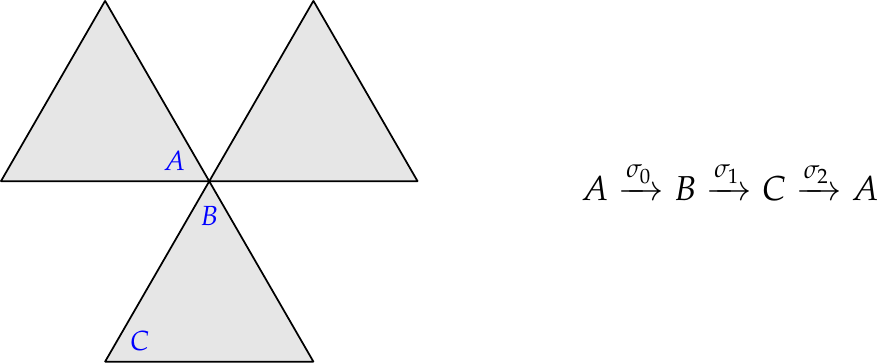}
\end{center}

\begin{example}
The following triples of permutations define connected 3-hypermaps of type $(0,2)$ and $(1,1)$, respectively.

\begin{align*}
\sigma_0 &= ( ~ 1 ~ ~ 12 ~ ) ~ ( ~ 2 ~ ~ 4 ~ ) ~ ( ~ 3 ~ ) ~ ( ~ 5 ~ ) ~ ( ~ 6 ~ ~ 7 ~ ) ~ ( ~ 8 ~ ) ~ ( ~ 9 ~ ~ 13 ~ ~ 10 ~ ) ~ ( ~ 12 ~ ) ~ ( ~ 14 ~ ) ~ ( ~ 15 ~ ) \\
\sigma_1 &= ( ~ 1 ~ ~ 2 ~ ~ 3 ~ ) ~ ( ~ 4 ~ ~ 5 ~ ~ 6 ~ ) ~ ( ~ 7 ~ ~ 8 ~ ~ 9 ~ ) ~ ( ~ 10 ~ ~ 11 ~ ~ 12 ~ ) ~ ( ~ 13 ~ ~ 14 ~ ~ 15 ~ ) \\
\sigma_2 &= ( ~ 1 ~ ~ 11 ~ ~ 10 ~ ~ 15 ~ ~ 14 ~ ~ 13 ~ ~ 8 ~ ~ 7 ~ ~ 5 ~ ~ 4 ~ ) ~ ( ~ 2 ~ ~ 6 ~ ~ 9 ~ ~ 12 ~ ~ 3 ~ ) \\
\\
\sigma_0 &= ( ~ 1 ~ ) ~ ( ~ 2 ~ ~ 5 ~ ) ~ ( ~ 3 ~ ~ 11 ~ ) ~ ( ~ 4 ~ ) ~ ( ~ 6 ~ ~ 7 ~ ) ~ ( ~ 8 ~ ~ 13 ~ ) ~ ( ~ 9 ~ ~ 10 ~ ) ~ ( ~ 12 ~ ~ 17 ~ ) ~ ( ~ 14 ~ ) ~ ( ~ 15 ~ ~ 16 ~ ) ~ ( ~ 18 ~ ) \\
\sigma_1 &= ( ~ 1 ~ ~ 2 ~ ~ 3 ~ ) ~ ( ~ 4 ~ ~ 5 ~ ~ 6 ~ ) ~ ( ~ 7 ~ ~ 8 ~ ~ 9 ~ ) ~ ( ~ 10 ~ ~ 11 ~ ~ 12 ~ ) ~ ( ~ 13 ~ ~ 14 ~ ~ 15 ~ ) ~ ( ~ 16 ~ ~ 17 ~ ~ 18 ~ ) \\
\sigma_2 &= ( ~ 1 ~ ~ 3 ~ ~ 10 ~ ~ 8 ~ ~ 15 ~ ~ 18 ~ ~ 17 ~ ~ 11 ~ ~ 2 ~ ~ 4 ~ ~ 6 ~ ~ 9 ~ ~ 12 ~ ~ 16 ~ ~ 14 ~ ~ 13 ~ ~ 7 ~ ~ 5 ~ )
\end{align*}
The figures below show their corresponding graphical representations. The hypermap below left is embedded in the sphere, while the hypermap below right is embedded in the torus obtained by gluing together the parallel edges of the rectangle.
\begin{figure}[ht!]
\centering
\includegraphics{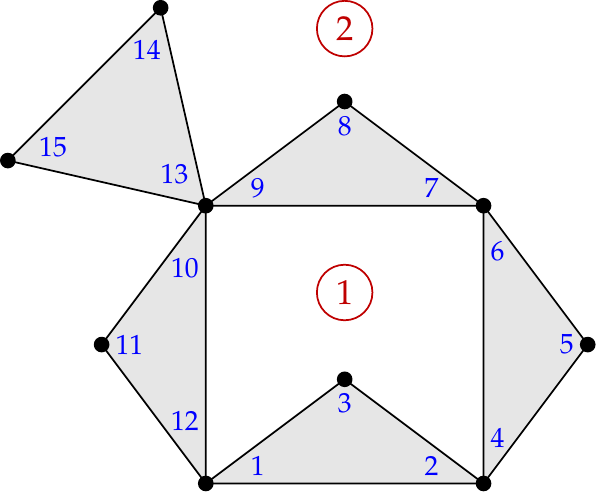}
\hspace{60pt}
\includegraphics{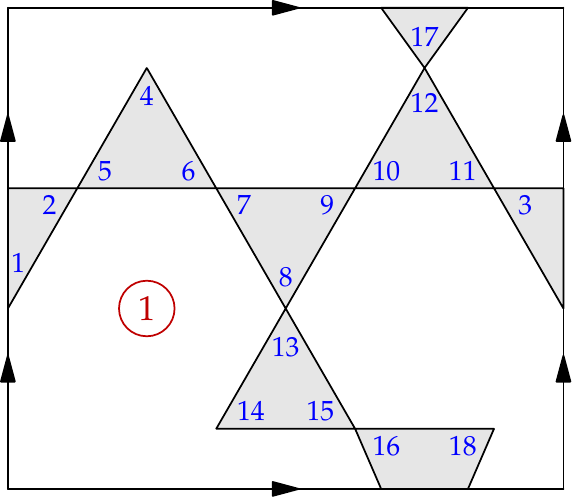}
\end{figure}

There are no nontrivial automorphisms of the first hypermap and one non-trivial automorphism of the second, corresponding to conjugation by the following permutation.
\[
( ~ 1 ~ ~ 4 ~ ) ~ ( ~ 2 ~ ~ 5 ~ ) ~ ( ~ 3 ~ ~ 6 ~ ) ~ ( ~ 7 ~ ~ 11 ~ ) ~ ( ~ 8 ~ ~ 12 ~ ) ~ ( ~ 9 ~ ~ 10 ~ ) ~ ( ~ 13 ~ ~ 17 ~ ) ~ ( ~ 14 ~ ~ 18 ~ ) ~ ( ~ 15 ~ ~ 16 ~ )
\]
Therefore, the first hypermap contributes 1 to the enumeration $M_{0,2}^{[3]}(5,10) = \frac{336}{5}$, while the second contributes $\frac{1}{2}$ to the enumeration $M_{1,1}^{[3]}(18) = \frac{52598}{3}$.
\end{example}

\section{Quantum curves for the enumeration of hypermaps} \label{main-proof}

\subsection{Statement of the main theorem}

We define the {\em free energies} for the enumeration of $a$-hypermaps to be the following generating functions.
\[
F_{g,n}^{[a]}(x_1, x_2, \ldots, x_n) = \sum_{b_1, b_2, \ldots, b_n = 1}^\infty M_{g,n}^{[a]}(b_1, b_2, \ldots, b_n) \, x_1^{-b_1} x_2^{-b_2} \cdots x_n^{-b_n}
\]
A subtle exception to the equation above is required in the case $(g,n) = (0,1)$.
\[
F_{0,1}^{[a]}(x_1) = -\log x_1 + \sum_{b_1 = 1}^\infty M_{0,1}^{[a]}(b_1) \, x_1^{-b_1}
\]
The {\em partition function} is defined from the free energies in the following way.
\[
Z^{[a]}(x, \h) = \exp \left[ \sum_{g=0}^\infty \sum_{n=1}^\infty \frac{\h^{2g-2+n}}{n!} \, F_{g,n}(x, x, \ldots, x) \right]
\]

\begin{remark}
Combinatorial justification for the definition of $F_{0,1}^{[a]}(x_1)$ stems from the fact that we have
\[
(-1)^n \, x_1 x_2 \cdots x_n \, \frac{\partial}{\partial x_1} \frac{\partial}{\partial x_2} \cdots \frac{\partial}{\partial x_n} F_{g,n}(x_1, x_2, \ldots, x_n) = \sum_{b_1, b_2, \ldots, b_n = 0}^\infty C_{g,n}(b_1, b_2, \ldots, b_n) \, x_1^{-b_1} x_2^{-b_2} \cdots x_n^{-b_n},
\]
where $C_{g,n}(b_1, b_2, \ldots, b_n)$ can be interpreted as the number of hypermaps with a distinguished choice of vertex in each face. In general, $C_{g,n}(b_1, b_2, \ldots, b_n) = b_1 b_2 \cdots b_n \, M_{g,n}(b_1, b_2, \ldots, b_n)$, although there is also the exceptional case $C_{0,1}(0) = 1$ corresponding to an isolated vertex on the sphere.
\end{remark}

In this section, we prove Theorem~\ref{th:qcurve}, which states that the partition function for the enumeration of $a$-hypermaps satisfies the following Schr\"{o}dinger-like equation, where $\widehat{x} = x$ and $\widehat{y} = -\h \frac{\partial}{\partial x}$.
\[
\left[ \widehat{y}^a - \widehat{x} \widehat{y} + 1 \right] Z^{[a]}(x, \h) = 0
\]

There is an issue of interpretation for this statement, since the definition of the partition function does not yield a convergent power series in $x$ and $\h$. In order to circumvent this issue, we define
\[
\overline{F}_{g,n}^{[a]}(x_1, x_2, \ldots, x_n) = \sum_{b_1, b_2, \ldots, b_n = 1}^\infty M_{g,n}^{[a]}(b_1, b_2, \ldots, b_n) \, x_1^{-b_1} x_2^{-b_2} \cdots x_n^{-b_n}
\]
for all $g \geq 0$ and $n \geq 1$. These are the usual free energies, without the exceptional logarithmic term that appears in the case $(g,n) = (0,1)$. This allows us to write the partition function in the following way.
\begin{equation} \label{eq:zbar}
Z^{[a]}(x, \h) = x^{-1/\h} \exp \left[ \sum_{g=0}^\infty \sum_{n=1}^\infty \frac{\h^{2g-2+n}}{n!} \, \overline{F}_{g,n}(x, x, \ldots, x) \right] = x^{-1/\h} \, \overline{Z}^{[a]}(x, \h)
\end{equation}
Note that the coefficient of $x^{-b}$ in $\overline{Z}(x, \h)$ is non-zero only for $b$ a positive integer, and is a Laurent polynomial in $\hbar$ --- in other words, $\overline{Z}(x, \h) \in \mathbb{Q}[\h^{\pm 1}][[x^{-1}]]$. Therefore, we may interpret Theorem~\ref{th:qcurve} in the following way.
\[
x^{1/\h} \left[ \widehat{y}^a - \widehat{x} \widehat{y} + 1 \right] x^{-1/\h} \, \overline{Z}(x, \h) = 0
\]
Explicitly applying the differential operator and simplifying yields the following rigorous interpretation of Theorem~\ref{th:qcurve}.

\begin{theorem} \label{th:mainmod}
The modified partition function $\overline{Z}^{[a]}(x, \h) \in \mathbb{Q}[\h^{\pm 1}][[x^{-1}]]$ of equation~(\ref{eq:zbar}) satisfies
\begin{equation} \label{eq:mainmod}
- \h x \frac{\partial}{\partial x} \overline{Z}^{[a]}(x, \h) = (-\h)^a \sum_{k=0}^a \binom{a}{k} ( -\h^{-1} )_k \, x^{-k} \left[ \frac{\partial^{a-k}}{\partial x^{a-k}} \overline{Z}^{[a]}(x, \h) \right].
\end{equation}
\end{theorem}

We use here the Pochhammer symbol to represent a falling factorial, which can be expressed in terms of unsigned Stirling numbers of the first kind.
\begin{equation} \label{eq:stirling}
(t)_k = t (t-1) (t-2) \cdots (t-k+1) = (-1)^k \sum_{j=0}^k \genfrac{[}{]}{0pt}{}{k}{j} (-t)^j
\end{equation}
Recall that the unsigned Stirling number of the first kind $\left[ \genfrac{}{}{0pt}{}{k}{j}  \right]$ is equal to the number of permutations in the symmetric group $S_k$ that have $j$ disjoint cycles. By convention, we set $\left[ \genfrac{}{}{0pt}{}{0}{0} \right] = 1$.

\subsection{Proof of the main theorem}

The logarithm of the modified partition function can be expressed in the following way.
\begin{align} \label{eq:logz}
\log \overline{Z}^{[a]}(x, \h) =& \sum_{g=0}^\infty \sum_{n=1}^\infty \frac{\h^{2g-2+n}}{n!} \, \overline{F}_{g,n}(x, x, \ldots, x) \nonumber \\
=& \sum_{g=0}^\infty \sum_{n=1}^\infty \frac{\h^{2g-2+n}}{n!} \, \sum_{b_1, b_2, \ldots, b_n = 1}^\infty M_{g,n}^{[a]}(b_1, b_2, \ldots, b_n) \, x^{-(b_1 + b_2 + \cdots + b_n)} \nonumber \\
=& \sum_{v=1}^\infty \sum_{e=1}^\infty f(v,e) \, \h^{(a-1)e - v} x^{-ae}
\end{align}
Here, $f(v, e)$ denotes the weighted count of connected, unlabelled $a$-hypermaps with $v$ vertices and $e$ hyperedges. As usual, the weight attached to a hypermap is equal to the reciprocal of its number of automorphisms. The last equality above uses the Riemann--Hurwitz calculation $v = 2-2g-n + (a-1) \frac{d}{a}$ that appears in Definition~\ref{def:hypermap} and the fact that the number of hyperedges is $\frac{d}{a}$. The factor of $\frac{1}{n!}$ disappears as we pass from the enumeration of labelled hypermaps to the enumeration of unlabelled hypermaps.

Although we have expressed the modified partition function $\overline{Z}^{[a]}(x, \h)$ in terms of the enumeration of connected hypermaps, it will be advantageous to express it in terms of the enumeration of possibly disconnected hypermaps. This motivates us to let $f^\bullet(v, e)$ denote the weighted count of possibly disconnected, unlabelled $a$-hypermaps with $v$ vertices and $e$ hyperedges. The particular form of the generating function above makes it amenable to the standard exponential trick to pass from a connected generating function to its possibly disconnected counterpart. For completeness, we provide the explicit verification below.

\begin{lemma} \label{lem:exponential}
The modified partition function $\overline{Z}^{[a]}(x, \h) \in \mathbb{Q}[\h^{\pm 1}][[x^{-1}]]$ of equation~(\ref{eq:zbar}) satisfies
\[
\overline{Z}^{[a]}(x, \h) = 1 + \sum_{e=1}^\infty \sum_{v=1}^\infty f^\bullet(v, e) \, \h^{(a-1)e - v}x^{-ae}.
\]
\end{lemma}

\begin{proof}
Exponentiate the formal power series of equation~(\ref{eq:logz}).
\begin{align*}
\overline{Z}^{[a]}(x, \h) &= \exp \left[ \sum_{e=1}^\infty \sum_{v=1}^\infty f(v, e) \, \h^{(a-1)e - v} \, x^{-ae} \right] \\
&= 1 + \sum_{e = 1}^\infty\sum_{v=1}^\infty \sum_{k=1}^\infty \frac{1}{k!}\h^{(a-1)e - v}x^{-ae} \mathop{\sum_{e_1 + e_2 + \cdots + e_k = e}}_{v_1 + v_2 + \cdots + v_k = v} \prod_{i=1}^k f(v_i, e_i)
\end{align*}
Here, we have used the fact that the number of vertices and hyperedges is additive over the disjoint union of connected hypermaps. The desired result is then equivalent to the following.
\[
f^\bullet(v, e) = \sum_{k=1}^\infty \frac{1}{k!} \mathop{\sum_{e_1 + e_2 + \cdots + e_k = e}}_{v_1 + v_2 + \cdots + v_k = v} \prod_{i=1}^k f(v_i, e_i).
\]
This equation reflects the fact that a possibly disconnected hypermap consists of $k$ connected components for some positive integer $k$. The factor of $\frac{1}{k!}$ kills the overcounting due to the fact that the connected components are unordered.
\end{proof}

We are now in a position to give a combinatorial proof of the main theorem of the paper.

\begin{proof}[Proof of Theorem~\ref{th:qcurve} and Theorem~\ref{th:mainmod}]
Use Lemma~\ref{lem:exponential} to express the left hand side of equation~(\ref{eq:mainmod}) as follows.
\begin{align*}
-x \h \frac{\partial}{\partial x} \overline{Z}(x, \h) &= - x \h \frac{\partial}{\partial x} \left[ 1 + \sum_{v=1}^\infty \sum_{e=1}^\infty f^\bullet(v,e) \, \h^{(a-1)e-v} \, x^{-ae} \right] \\
&= \sum_{v=1}^\infty \sum_{e=1}^\infty ae f^\bullet(v,e) \h^{(a-1)e-v+1} x^{-ae}
\end{align*}

Now use Lemma~\ref{lem:exponential} to express the right hand side of equation~(\ref{eq:mainmod}) as follows.
\begin{align*}
&~ (-\h)^a \sum_{k=0}^a \binom{a}{k} ( -\h^{-1})_{k} \, x^{-k} \left[ \frac{\partial^{a-k}}{\partial x^{a-k}} \overline{Z}^{[a]}(x, \h) \right] \\
=&~ (-\h)^a (-\h^{-1})_a \, x^{-a} + (-\h)^a \sum_{k=0}^a \binom{a}{k} (-\h^{-1})_k \, x^{-k} \left[ \frac{\partial^{a-k}}{\partial x^{a-k}} \left( \sum_{v=1}^\infty \sum_{e=1}^\infty f^\bullet(v,e) \, \h^{(a-1)e-v} \, x^{-ae} \right) \right] \\
=&~ (-\h)^a (-\h^{-1})_a \, x^{-a} + \sum_{k=0}^a \binom{a}{k} \sum_{j=0}^k \left[ \genfrac{}{}{0pt}{}{k}{j} \right] \left[ \sum_{v=1}^\infty \sum_{e=1}^\infty f^\bullet(v,e) \, \h^{(a-1)e-v+a-j} \, x^{-ae-a} \frac{(ae+a-k-1)!}{(ae-1)!} \right] \\
=&~ x^{-a} \sum_{j=0}^a \left[ \genfrac{}{}{0pt}{}{a}{j} \right] \h^{a-j} + \sum_{k=0}^a \binom{a}{k} \sum_{j=0}^k \left[ \genfrac{}{}{0pt}{}{k}{j} \right] \left[ \sum_{v=1+j}^\infty \sum_{e=2}^\infty f^\bullet(v-j,e-1) \, \h^{(a-1)e-v+1} \, x^{-ae} \frac{(ae-k-1)!}{(ae-a-1)!} \right]
\end{align*}

We can use the fact that $f^\bullet(v-j,e-1) = 0$ for $v-j \leq 0$ and $e \geq 2$ to write this in the following way. The first term above corresponds precisely to the $e = 1$ summand below, as long as we take $f^\bullet(0,0) = 1$ as the only non-zero value of $f^\bullet(v, e)$ with $v=0$ or $e=0$. We also change the order of summation for convenience.
\[
\sum_{v=1}^\infty \sum_{e=1}^\infty \sum_{k=0}^a \sum_{j=0}^k \binom{a}{k} \left[ \genfrac{}{}{0pt}{}{k}{j} \right] \frac{(ae-k-1)!}{(ae-a-1)!} f^\bullet(v-j,e-1) \, \h^{(a-1)e-v+1} \, x^{-ae}
\]

Now we may equate the coefficients on both sides of equation~(\ref{eq:mainmod}) to obtain the following.
\[
ae f^\bullet(v,e) = \sum_{k=0}^a \sum_{j=0}^k \binom{a}{k} \left[ \genfrac{}{}{0pt}{}{k}{j} \right] \frac{(ae-k-1)!}{(ae-a-1)!} \, f^\bullet(v-j,e-1)
\]

In order to prove this combinatorially, we define $F^\bullet(v,e) = f^\bullet(v,e) \times (ae)!$ to be the number of possibly disconnected, unlabelled $a$-hypermaps with $v$ vertices and $e$ hyperedges, with a total ordering on the angles. One can think of the angles as being labelled from 1 up to $ae$. So it now suffices to prove the following equation.
\begin{equation} \label{eq:recursion}
F^\bullet(v,e) = \sum_{k=0}^a \sum_{j=0}^k \frac{(ae-1)!}{(ae-a)!} \binom{a}{k} \left[ \genfrac{}{}{0pt}{}{k}{j} \right] \frac{(ae-k-1)!}{(ae-a-1)!} \, F^\bullet(v-j,e-1)
\end{equation}

By definition, the left hand side of equation~(\ref{eq:recursion}) counts angle-ordered, possibly disconnected, unlabelled $a$-hypermaps with $v$ vertices and $e$ hyperedges. We call the hyperedge containing the largest angle label {\em marked}. We will show that the $(k, j)$ summand on the right hand side of equation~(\ref{eq:recursion}) counts such hypermaps where
\begin{itemize}
\item the marked hyperedge has $k$ angles that are not adjacent to any other hyperedge; and
\item these $k$ angles are glued together to form $j$ vertices.
\end{itemize}
Clearly, for any angle-ordered, unlabelled $a$-hypermap, there exists a unique pair $(k, j)$ for which the two conditions above hold. We remark that the $(0, 0)$ summand does indeed contribute to the right hand side, since $\left[ \genfrac{}{}{0pt}{}{0}{0} \right] = 1$.

In order to obtain equation~(\ref{eq:recursion}), consider removing the marked hyperedge.
\begin{itemize}
\item The factor $F^\bullet(v-j,e-1)$ is the number of hypermaps that can remain once the marked hyperedge is removed.
\item The factor $\frac{(ae-1)!}{(ae-a)!}$ is the number of ways to choose the remaining angle labels on the marked hyperedge.
\item The factor $\binom{a}{k}$ is the number of ways to choose the $k$ angles on the marked hyperedge that are not adjacent to any other hyperedge.
\item The factor $\left[ \genfrac{}{}{0pt}{}{k}{j} \right]$ is the number of ways that these $k$ angles can be glued together to form $j$ vertices. To see this, recall that $\left[ \genfrac{}{}{0pt}{}{k}{j} \right]$ is the number of permutations in $S_k$ that have $j$ disjoint cycles. Each cycle determines a set of angles that are glued together to form a vertex, as well as the cyclic orientation of the angles at that vertex.
\item The factor $\frac{(ae-k-1)!}{(ae-a-1)!}$ is the number of ways that the remaining $a-k$ vertices on the marked hyperedge can connect to the rest of the hypermap. To see this, note that there are $ae-a$ places to glue the first of these angles, $ae-a+1$ places to glue the second, $ae-a+2$ places to glue the third, and so on.
\end{itemize}
Therefore, we have shown that both sides of equation~(\ref{eq:recursion}) are equal to the number of angle-ordered, possibly disconnected, unlabelled $a$-hypermaps with $v$ vertices and $e$ hyperedges.
\end{proof}

As mentioned earlier, the $a=2$ case of Theorem~\ref{th:mainmod} recovers the quantum curve for the enumeration of ribbon graphs, which was rigorously determined by Mulase and Su{\l}kowski~\cite{mul-sul}. We emphasise that our proof bypasses the calculations in their paper by interpreting the partition function as a generating function for disconnected unlabelled objects. This viewpoint may prove valuable for quantum curves associated to other enumerative problems. One advantage of our proof is that it yields a concise explicit formula for the partition function in terms of falling factorials.

\begin{proposition}
The modified partition function $\overline{Z}^{[a]}(x, \h) \in \mathbb{Q}[\h^{\pm 1}][[x^{-1}]]$ of equation~(\ref{eq:zbar}) satisfies
\[
\overline{Z}^{[a]}(x, \h) = \sum_{e=0}^\infty \frac{( -\h^{-1} )_{ae}}{e!} \left[ (-1)^a \frac{\h^{a-1}}{ax^a} \right]^e.
\]
\end{proposition}

\begin{proof}
Note that $F^\bullet(v, e)$ is equal to the number of triples of permutations $(\sigma_0, \sigma_1, \sigma_2)$ in the symmetric group $S_{ae}$ such that $\sigma_0$ has $v$ cycles, $\sigma_1$ has cycle type $(a, a, \ldots, a)$ with $e$ cycles, and $\sigma_0 \sigma_1 \sigma_2 = \text{id}$. The permutation $\sigma_0$ can be chosen in $\left[ \genfrac{}{}{0pt}{}{ae}{v} \right]$ ways, while the permutation $\sigma_1$ can be chosen in $\frac{(ae)!}{a^e \times e!}$ ways. Once these have been selected, the permutation $\sigma_2 = (\sigma_0 \sigma_1)^{-1}$ is uniquely defined. Therefore, we have
\[
f^\bullet(v, e) = \frac{1}{a^e \times e!} \left[ \genfrac{}{}{0pt}{}{ae}{v} \right],
\]
and combining this with Lemma~\ref{lem:exponential} and equation~(\ref{eq:stirling}) produces the desired  expression.
\end{proof}

\begin{remark}
It is a curious fact that $\overline{Z}^{[a]}(x, -1) = 1$, which is equivalent to the following statement for a fixed positive integer $e$.
\[
\sum_{v \text{ even}} f^\bullet(v, e) = \sum_{v \text{ odd}} f^\bullet(v, e)
\]
This is a consequence of the simple fact that the function sending the hypermap $(\sigma_0, \sigma_1, \sigma_2)$ to the hypermap $((1\,2) \circ \sigma_0, \sigma_1, \sigma_2 \circ (1\,2))$ is an involution on the set of hypermaps with $e$ hyperedges that either increases or decreases the number of vertices by 1.
\end{remark}

\section{A conjecture on hypermaps and topological recursion} \label{conjecture}

\subsection{Statement of the conjecture}

The semi-classical limit of the quantum curve $\widehat{y}^a - \widehat{x} \widehat{y} + 1$ is the spectral curve $y^a - xy + 1 = 0$. Introducing the rational parameter $z$ allows us to express this parametrically as $x(z) = z^{a-1} + \frac{1}{z}$ and $y(z) = z$. The theory of quantum curves and their relation to the topological recursion proposed by Gukov and Su{\l}kowski immediately suggest the following, which we introduced earlier as Conjecture~\ref{con:hypermaps}. We use here the notation $x_i = x(z_i)$ for $i = 1, 2, \ldots, n$.

\begin{conjecture*}
For a fixed positive integer $a$, the topological recursion applied to the rational spectral curve $x(z) = z^{a-1} + \frac{1}{z}$ and $y(z) = z$ produces correlation differentials whose expansions at $x_i = \infty$ satisfy
\[
\omega_{g,n} = \sum_{b_1, b_2, \ldots, b_n = 1}^\infty M_{g,n}^{[a]}(b_1, b_2, \ldots, b_n) \, \prod_{i=1}^n \frac{b_i}{x_i^{b_i+1}} \, \dd x_i \qquad \qquad \text{for } 2g-2+n>0.
\]
\end{conjecture*}

This provides one of few known instances in which the semi-classical limit of the quantum curve gives rise to a conjecture relating an enumerative problem with the topological recursion. Another instance of this phenomenon is the case of spin Hurwitz numbers, for which the quantum curve has been calculated, yet the connection to the topological recursion is still only conjectural~\cite{mul-sha-spi}.

There is considerable numerical evidence to support the above conjecture. For example, in the cases $(g,n) = (0,3), (0,4), (1,1), (1,2), (2,1)$ for $a = 1, 2, 3, 4, 5$, the values of $M_{g,n}^{[a]}(b_1, b_2, \ldots, b_n)$ are consistent with the correlation differential $\omega_{g,n}$ up to order 15. These calculations also support the following conjecture.

\begin{conjecture} \label{con:polynomial}
There exists a quasi-polynomial $P_{g,n}^{[a]}$ modulo $a$ of degree $3g-3+n$ such that
\[
M_{g,n}^{[a]}(b_1, b_2, \ldots, b_n) = \prod_{i=1}^n \binom{b_i-1}{\lfloor \frac{b_i-1}{a} \rfloor} \times P_{g,n}^{[a]}(b_1, b_2, \ldots, b_n).
\]
\end{conjecture}

To say that $P_{g,n}^{[a]}(b_1, b_2, \ldots, b_n)$ is a quasi-polynomial modulo $a$ means that it is a polynomial when restricted to each coset of $(a\mathbb{Z})^n \subseteq \mathbb{Z}^n$. This quasi-polynomial behaviour for the enumeration of hypermaps is analogous to the polynomial behaviour for simple Hurwitz numbers, which was first conjectured by Goulden and Jackson~\cite{gou-jac}. The first proof of this fact relied upon the ELSV formula, which expresses simple Hurwitz numbers in terms of the intersection theory of moduli spaces of curves~\cite{eke-lan-sha-vai}. A recent proof that circumvents such involved algebraic geometric considerations was recently found, relying instead on the technology of the semi-infinite wedge space~\cite{dun-kaz-ora-sha-spi}. A proof of Conjecture~\ref{con:polynomial} may help to shed further light on the polynomial behaviour of Hurwitz numbers and their generalisations.

\subsection{Evidence and applications}

The topological recursion produces all of the correlation differentials $\omega_{g,n}$ from the base cases $\omega_{0,1}$ and $\omega_{0,2}$. This viewpoint was enunciated by Dumitrescu et al., where they state: {\em the Laplace transform of the unstable geometries $(g,n) = (0,1)$ and $(0,2)$ determines the spectral curve}~\cite{dum-mul-saf-sor}. Since $\omega_{0,2}$ is canonical in the case of a genus zero spectral curve, one can expect to deduce the spectral curve for certain enumerative problems from the $(g,n) = (0,1)$ information alone. The following proposition verifies that Conjecture~\ref{con:hypermaps} is consistent with this statement.

\begin{proposition} \label{th:01-formula}
If $y^a - xy + 1 = 0$, then the expansion of $y$ at $x = \infty$ is given by the formula
\[
y = 1 + \sum_{b=1}^\infty b M_{0,1}^{[a]}(b) \, x^{-b-1}.
\]
\end{proposition}

\begin{proof}
Define $N_0 = 1$ and $N_b = b M_{0,1}^{[a]}(b)$ for all positive integers $b$. We interpret $N_b$ combinatorially as the number of connected $a$-hypermaps of type $(0,1)$ with one marked angle and whose face has perimeter $b$. Removing the hyperedge with the marked angle leaves an $a$-tuple of hypermaps of type $(0,1)$, each marked by the unique angle adjacent to the removed hyperedge. Therefore, we have deduced the following recursion, which can be used to calculate all values of $N_b$ from the base cases $N_0 = 1$ and $N_1 = N_2 = \cdots = N_{a-1} = 0$.
\[
N_b = \sum_{m_1+m_2+\cdots+m_a=b-a} N_{m_1} N_{m_2} \cdots N_{m_a}
\]

Now suppose that $Y = \displaystyle\sum_{b=0}^\infty N_b \, X^{-b-1}$ and use the recursion above to conclude that
\[
Y^a = \sum_{m_1, m_2, \ldots, m_a=0}^\infty N_{m_1} N_{m_2} \cdots N_{m_a} \, X^{-(m_1+m_2+\cdots+m_a+a)} = \sum_{m=0}^\infty N_{m+a} \, X^{-m-a} = \sum_{m=a}^\infty N_m \, X^{-m} = XY - 1.
\]
Therefore, $Y^a - XY + 1 = 0$ and the desired result follows immediately.
\end{proof}

\begin{remark}
In fact, one can use Proposition~\ref{th:01-formula} to recover the following formula.
\[
bM_{0,1}^{[a]}(b) = \begin{cases} \frac{a}{ab + a - b} \binom{b}{b/a}, & \text{if } b \equiv 0 \pmod{a} \\
0, & \text{if } b \not\equiv 0 \pmod{a}. \end{cases}
\]
These values can be considered as generalisations of the Catalan numbers, which are recovered in the case $a = 2$, where we have $C_m = 2m M_{0,1}^{[2]}(2m)$.
\end{remark}

The correlation differentials produced by the topological recursion satisfy {\em string equations}~\cite{eyn-ora07a}.
\begin{eqnarray*}
\sum_\alpha \mathop{\text{Res}}_{z=\alpha} y(z) \, \omega_{g,n+1}(z, z_S) &=& - \sum_{i=1}^n \dd z_i \, \frac{\partial}{\partial z_i} \left( \frac{\omega_{g,n}(z_S)}{\dd x(z_i)} \right) \\
\sum_\alpha \mathop{\text{Res}}_{z=\alpha} x(z) y(z) \, \omega_{g,n+1}(z, z_S) &=& - \sum_{i=1}^n \dd z_i \, \frac{\partial}{\partial z_i} \left( \frac{x(z_i) \omega_{g,n}(z_S)}{\dd x(z_i)} \right)
\end{eqnarray*}
They are also known to satisfy the following {\em dilaton equation}, where $\Phi(z) = \int y(z) \, \dd x(z)$.
\[
\sum_\alpha \mathop{\text{Res}}_{z=\alpha} \Phi(z) \, \omega_{g,n+1}(z, z_S) = (2g-2+n) \, \omega_{g,n}(z_S)
\]

In conjunction with Conjecture~\ref{con:hypermaps}, these equations should yield non-trivial relations between values of $M_{g,n+1}^{[a]}(b_1, b_2, \ldots, b_{n+1})$ and values of $M_{g,n}^{[a]}(b_1, b_2, \ldots, b_n)$. For the spectral curve $x(z) = z + \frac{1}{z}$ and $y(z) = z$, the dilaton equation was used to define the enumeration of ribbon graphs where some of the $b_i$ are equal to zero~\cite{nor13}. This was a crucial ingredient in the enumeration of lattice points in the Deligne--Mumford compactification of moduli spaces of curves~\cite{do-nor}.

The dilaton equation furthermore allows one to define the symplectic invariants $F_g = \omega_{g,0} \in \mathbb{C}$ for $g \geq 2$. It was shown by Norbury that the symplectic invariants for the spectral curve $x(z) = z + \frac{1}{z}$ and $y(z) = z$ satisfy $F_g = \chi({\mathcal M}_g)$, the orbifold Euler characteristic of the moduli space of genus $g$ curves~\cite{nor13}. It would be interesting to calculate the symplectic invariants for the sequence of spectral curves introduced in Conjecture~\ref{con:hypermaps} and to determine whether there is an analogous connection to the geometry of moduli spaces.

\begin{small}
\bibliographystyle{hacm}
\bibliography{hypermaps-qcurve}
\nocite{*}

\textsc{School of Mathematical Sciences, Monash University, VIC 3800, Australia} \\
{\em Email:} \href{mailto:norm.do@monash.edu}{norm.do@monash.edu}

\textsc{School of Mathematics and Statistics, The University of Sydney, NSW 2006, Australia} \\
{\em Email:} \href{mailto:david.manescu@gmail.com}{david.manescu@gmail.com}
\end{small}

\end{document}